\documentclass[preprint,12pt,authoryear]{elsarticle}

\usepackage{amssymb}
\usepackage{amsmath}
\newtheorem{thm}{Theorem}
 \newtheorem{lem}[thm]{Lemma}
 \newdefinition{rmk}{Remark}
 \newproof{pf}{Proof}
\journal{Applied Mathematics Letters}

\begin{document}

\begin{frontmatter}


 %\author{Anshui Li\corref{cor1}\fnref{label2}}
  \author{Anshui Li}
 \ead{anshuili@hznu.edu.cn}
 \author{Yong Chen\corref{cor1}}
 \ead{ychen@hznu.edu.cn}
 %\fntext[label2]{anshuili@hznu.edu.cn}
 \cortext[cor1]{Corresponding author}
% \fntext[label3]{Correspondence anthor}

%\title{Convergence of Kronecker Graph Processes via Differential Equations with Unique Solutions}
%\title{Convergence of P{\'o}lya urn  model via Wormald's equations techniques}
\title{Convergence of  Coupon Collecting Process via Wormald's  Differential Equation Method}
 %\author{Anshui Li, Yong Chen}
 \address{Department of Mathematics, School of Science, Hangzhou Normal University, Hangzhou, Zhejiang, P.R.China}

%\author{Anshui Li, Yong Chen}

%\address{Department of Mathematics, School of Science, Hangzhou Normal University, Hangzhou, Zhejiang, P.R.China}

\begin{abstract}
%One novel epidemic model called SIR with the influence of Media is studied in this paper. Compare the models in the literature, ours is much more feasible.
To approximate the trajectories of a stochastic process by the solution of some differential equations is widely used in the fields of probability, computer science and combinatorics. In this paper, the convergence of coupon collecting Process is studied via the differential equation techniques originally proposed in \cite{wormald1995differential}  and modified in \cite{warnke2019wormald}. In other words, we give a novel approach to analysis the classical coupon collector's problem and keep track of all coupons in the process of collecting. .
\end{abstract}



\begin{keyword}
Wormald's  Equation Method \sep Coupon Collector's Problem \sep Differential Equations\sep Stochastic Processes\sep P{\'o}lya Urn Model

 \MSC[2010] 34F05\sep 34Exx \sep 60H10 
 %code 
 %(2000 is the default)

\end{keyword}

\end{frontmatter}


\section{Introduction}
\label{}
%The study of epidemics has a fairly long history, with a great number of successful mathematical models available in the literature. Many techniques are introduced to study epidemic models, especially the dynamical systems,

The differential equations normally arise from the dynamic of the stochastic processes concerned. It is very natural and powerful to approximate the trajectories of a stochastic process by the solutions of differential equations. This technique dates back to the pioneering work in the field of applied probability for continuous-time Markov process in \cite{kurtz1970solutions,kurtz1981approximation}; and later were popularized in the combinatorics, especially in the study of evolution of random structures in~\cite{karp1981maximum,wormald1999differential,riordan2016convergence,warnke2016method,heredia2018modelling}.

The coupon collector's problem is one of very classical models in applied probability and combinatorics, see~\cite{glavavs2018new, boneh1997coupon} for some general information. In this paper, we will focus on one novel approach of differential equation to consider the classical coupon collector's problem and get some interesting result.
The paper is organized as follows: the differential equation techniques for random structures is reviewed briefly in Section 2. Our results on the convergence of coupon collecting process will be showed and proved in Section 3. In the last Section, some remarks and possible extensions are given.

\section{The differential equation technique for random structures}
In this section, we will briefly review the differential equation techniques given by \cite{wormald1995differential, wormald1999differential} and is simply modified and proved again in \cite{warnke2019wormald}. Roughly speaking, this technique shows some general criteria to guarantee that some given parameters in a family of discrete random processes converge to the solution of a system of differential equations. 

We will follow the notations from the original paper in \cite{wormald1995differential}. This technique is designed for discrete time stochastic processes(mostly for random graph processes) at the very beginning. Such a process is a probability space $\Omega$ which is denoted by $(Q_0,~Q_1,\cdots)$, where each $Q_i$ takes values in some set $S$. The elements of $\Omega$ are sequences $(q_0,q_1,\cdots )$, where each $q_i\in S$. We use $H_t$ to denote $(Q_0,Q_1,\cdots, Q_t)$, the history(filtration) of the process up to time $t$. For a function $y$ defined on histories, the random variable $y(H_t)$ is denoted by $Y_t$ for convenience.

Consider a sequence $\Omega_n,~n=1,2,...$ of random processes. Thus $q_t=q_t(n)$ and $S=S_n$, but the dependence on $n$ is usually dropped from the notation for convenience. An event occurs \emph{almost surely} if its probability in $\Omega_n$ is $1-o(1)$. 

The set of all $h_t=(q_0,q_1,\cdots,q_t)$ will be denoted by $S_n^+$, where each $q_i\in S_n,t=0,1,~\cdots$.

A function $f(u_1,\cdots,u_j)$ is called satisfies a \emph{Lipschitz condition} on $D\subset \mathbb{R}^j$ if for some constant $L>0$ exists with the following property holds
\[
|f(u_1,\cdots,u_j)-f(v_1,\cdots,v_j)|\le L\sum_{i=1}^j|u_i-v_i|
\] for all $(u_1,\cdots, u_j)$ and $(v_1,\cdots,v_n)\in D$.

It is time to state the main technique originally proposed in \cite{wormald1995differential} now. The theorem goes as follows:
\begin{thm}\label{thm:main}
Let a be fixed. For $1\le l \le a$, let $y^{(l)}:~\cup_n S_n^+\to \mathbb{R}$ and $f_l: \mathbb{R}^{a+1}\to \mathbb{R}$, such that for some constant $C$ and all $l$, $|y^{(l)}(h_t)|< Cn$ for all $h_t\in S_n^+$ for all $n$. Suppose also that for some function $m=m(n)$:
\begin{itemize}
\item there is a constant $C^{'} $ such that, for all $ t<m$ and all $l$,
\[
|Y_{t+1}^{(l)}-Y_{t}^{(l)}|< C^{'}
\] always;
\item for all $l$ and uniformly over all $t < m$,
\[
\mathbb{E}(Y_{t+1}^{(l)}-Y_t^{(l)}|H_t)=f_l(t/n,Y_t^{(1)}/n,\cdots,Y_t^{(a)}/n)+o(1)
\]always;
\item for each $l$ the function $f_l$ is continuous and satisfies a Lipschitz condition on $D$, where $D$ is some bounded connected open set containing the intersection of $\{(t,z^{(1)},\cdots,z^{(a)}):t\ge 0\}$ with some neighborhood of $\{(0,z^{(1)},\cdots,z^{(a)}): \mathbb{P}(Y_0^{(l)}=z^{(l)}n, 1\le l\le a)\neq 0~ for~some~n \}$
\end{itemize}
Then:
\begin{enumerate}
\item For $(0,\hat{z}^{(1)},\cdots,\hat{z}^{(a)})\in D$ the system of differential equations
\[
\frac{d z_l}{ds}=f_l(s,z_1,\cdots,z_a),~~~l=1,\cdots,a,
\]has a unique solution in $D$ for $z_l:~\mathbb{R}\to \mathbb{R}$ passing through 
\[
z_l(0)=\hat{z}^{(l)}, ~~1\le l\le a,
\]and which extends to points arbitrarily close to the boundary of $D$.
\item Almost surely
\[
Y_t^{(l)}=nz_l(t/n)+o(n)
\]uniformly for $0\le t \le\min\{\sigma n,m\}$ and for each $l$, where $z_l(t)$ is the solution in $(1)$ with $z^{(l)}={Y_0^{(l)}}/{n}$, and $\sigma=\sigma(n)$ is the supremum of those $s$ to which the solution can be extended.
\end{enumerate}
\end{thm}
\begin{rmk}
The theorem above can be used to determine various  quantities asymptotically almost surely in proper random graph processes and other related random processes. To some extent, the key point in the theorem above is to show the condition expectations can be expressed in a series of proper functions with "nice" properties.
\end{rmk}

\section{Coupon collector's problem revised}

In this section ,we will revise the coupon collector's problem via Wormald's equation techniques, which turns out to very easy to solve and can  be extended to other urn models. 

The coupon collector's problem statement goes as follows: there are $n$ different coupons needed to be collected. The coupons are obtained at random, drawing from a box that contains all coupons which are labelled from $1$ to $n$. When a coupon is collected the player keeps it, and the target is to collect at least one coupon of every type. There are many results related to this model from different facts. This model can find many powerful applications in many areas, see~\cite{zheng2019authorship,tan2019bitcoin,glavavs2018new,lundow2019revisiting}. There are many results of this classical problem from different point of views and a number of applications in analysis of algorithms and combinatorics can be found in \cite{mitzenmacher2017probability}.

For some technical reason, we will consider the coupon collector's problem as a stochastic process called coupon collecting process. Let $Y_t^i$ be the number of coupons with $i$ copys after $t$ steps, and $Y_t$ be the coupon collected at step $t$, and $T$ denote the first time that all the coupons are collected.
Then we have the following lemma which is mentioned in \cite{mitzenmacher2017probability}. 
\begin{lem}For any constant $c$, we have
\[
\lim_{n \to \infty}\mathbb{P}(T \ge n\ln n+cn)=1-e^{-e^{-e^c}}.
\]
\end{lem}
%We can choose $a=\ln n$ and $m=n\ln n+cn$ for $c>1$.  
\begin{rmk}
The lemma above states that, for large $n$, the number of coupons required should be close to $n\ln n$. That means, we can choose $m$ close to $n\ln $ for the coupon collecting process when one wants to use Theorem~\ref{thm:main}.
\end{rmk}
We obtain the following theorem for coupon collecting process.\begin{thm}
Almost surely, we have
\[
Y_t^i =z_i(t/n)n+o(n)
\]uniformly for $0\le t\le n\ln n$ and for each $0<i<l$ with $l$ some positive constant, where the $z_i(t)$ form the solution to
\[
\frac{dz_i}{ds}=f_i(s,z_0,z_1,\cdots,z_l),
\] in which $f_i(s,z_0,\cdots,z_l)={z_{i-1}-z_i}$, with initial conditions $z_i(0)=0$ for $i>0$ and $z_0=1$.
\end{thm}

\begin{pf} Since $Y_{t+1}^i$ and $Y_t^i$ is the number of coupons with $i+1$ and $i$ respectively up to time $t$. $Y_{t+1}^i$ change at most by $1$ comparing to $Y_{t}^i$ for every $i$, which means
\[
|Y_{t+1}^i-Y_t^i|\le 1.
\]
Then we calculate the conditional expectation for the difference of $Y_{t+1}^i-Y_t^i$ based on the filtration $H_t$. This variable will be $1$ if the coupon collected at time $t+1$ is one type with $i-1$ copys  at time $t$. The probability of this case is $y_{t}^{i-1}/n$. It is easy to see that it is $-1$ with probability $y_t^i/n$ with similar argument. To sum up, we have 
\[
\mathbb{E}(Y_{t+1}^i-Y_t^i|H_t)=\frac{y_t^{i-1}}{n}-\frac{y_t^i}{n}.
\]
And for each $i$, the function $f_i(s,z_0,\cdots,z_l)=z_{i-1}-z_i$ is obviously continuous and
\begin{align}
\notag
|f_i(u_0,u_1,\cdots,u_n)-f_i(v_0,v_1,\cdots,v_n)| &=|u_{i-1}-u_i-v_{i-1}+v_i|\\
\notag
&=|u_i-v_i|+|u_{i-1}-v_{i-1}|\\
&\le\sum_{i=0}^l|u_i-v_i|
\end{align}
i.e., function $f_i$ satisfies a \emph{Lipschitz condition} on some bounded connection open set.

The result above is obtained easily as we have checked the conditions in Theorem \ref{thm:main}, which concludes the proof.
\end{pf}
\begin{rmk}
The result above show that the coupons with same copys is uniformly determined by a series of differential equations.
 Roughly speaking, the theorem above states that all the information of coupons collected in the collecting process at any time $t$ can be tracked in a whole.
\end{rmk}
\section{Conclusion}
Many quantities related to random graph processes are approximately determined almost surely as the total number of vertices or time increases with the help of Wormald's differential equation method. In particular, most of the P{\'o}lya Urn model can be regarded as some graph process or some variants of coupon collector's problem, which turns out can be explored with Wormald's differential equation method  with the same procedure mentioned in this paper. We will try to explore more extensions which can be solved with Wormald's differential equations method in the follow-up work.



\section*{Acknowledgments}
This work is supported by the National Natural Science Foundation of China under Grant No.11901145 and Zhejiang Provincial Natural Science Foundation of China under Grant No.LQ18A010007.
\section*{References}
\bibliographystyle{elsarticle-harv}
 \bibliography{diff}

\end{document}